\documentclass[12pt,a4paper]{article}
\usepackage{amsfonts}

\usepackage{amssymb}

\usepackage{amsmath}

\newtheorem{theorem}{Theorem}[section]
\newtheorem{definition}[theorem]{Definition}
\newtheorem{lemma}[theorem]{Lemma}

\newtheorem{example}[theorem]{Example}

\begin{document}
\title{Gr\"{o}bner-Shirshov bases for some one-relator groups\footnote{Supported by the
NNSF of China (No.10771077) and the NSF of Guangdong Province
(No.06025062).} }
\author{
Yuqun Chen and Chanyan Zhong  \\
\\
{\small \ School of Mathematical Sciences}\\
{\small \ South China Normal University}\\
{\small \ Guangzhou 510631}\\
{\small \ P. R. China}\\
{\small \ yqchen@scnu.edu.cn} \\
{\small \ chanyanzhong@yahoo.com.cn}}
\date{}
\maketitle \noindent\textbf{Abstract:} In this paper, we  prove that
 two-generator one-relator groups with depth less than or equal to 3 can be effectively embedded into a tower
of HNN-extensions in which each group has the effective standard
normal form. We give an example to show how to deal with some
general cases for one-relator groups.  By using the Magnus' method
and Composition-Diamond lemma, we reprove the G. Higman, B. H.
Neumann and H. Neumann's embedding theorem.

\noindent \textbf{Key words: } Group, HNN-extension,
Gr\"{o}bner-Shirshov basis, Standard normal form.

\noindent {\bf AMS} Mathematics Subject Classification(2000): 20E06,
20F05, 16S15, 13P10

\section{Introduction}

Higman, Neumann, Neumann (1949, \cite{hnn}, see also R. C. Lyndon
and P. E. Schupp \cite{ly}, p.188) proved that any countable group
with $\leq n$ relations can be effectively embedded into a group
generated by two elements with $\leq n$ relations. Even before, this
result was known for one-relator groups (W. Magnus \cite{mu31}, see
also W. Magnus, A. Karrass and D. Solitar \cite{mks}, p.259). So,
any finitely generated one-relator group can be effectively embedded
into one-relator group with two generators,

$$
G=gp\langle x,\ y|x^{n_1}y^{m_1}\cdots x^{n_k}y^{m_k}=1\rangle,
$$
where $n_i, \ m_i\neq0, \ k\geq0$. We call  $k$  the depth of $G$.

On the other hand, any one-relator group can be effectively embedded
into a tower of HNN-extensions essentially by the Magnus' method
(see \cite{ly}, p.198, and Moldavanskii \cite{mo67}). There is a
conjecture, stated by L. A. Bokut, that each group of the Magnus
tower for any one-relator group is a group with the standard normal
form in the sense of Bokut \cite{bo66} (see also \cite{b94}). If it
is true, it would give another proof for the decidability of the
word problem for any one-relator group.

In this paper, we  prove that one-relator groups with depth less
than or equal to 3 can be effectively embedded into a tower of
HNN-extensions in which each group has the effective standard normal
form. We give an example to show how to deal with some general cases
for one-relator groups. We  use the Magnus' method and
Composition-Diamond lemma to reprove the G. Higman, B. H. Neumann
and H. Neumann's embedding theorem.

\section{Standard normal form and Gr\"obner-Shirshov basis}
In this section, we will cite some literatures about the definition
of groups with the standard normal form and Composition-Diamond
lemma on free associative algebra $k\langle X\rangle$.

\begin{definition}\label{d2.1}\cite{ly} \
Let $G$ be a group, $A,\ B$ the subgroups of $G$ with $\phi:
A\rightarrow B$ an isomorphism. Let
$$
H=gp\langle G,t|\ t^{-1}at=b,\ a\in A,\ b=\phi(a)\rangle.
$$
Then $H$ is called an HNN-extension of $G$ relative to $A,\ B \mbox{
and } \phi$.
\end{definition}

\begin{definition}\label{d2.1}\cite{n52, n54, n55} \
Let $G$ be a group, $t$ a letter, $A_i,\ B_i\in G,\ \ i\in I$. Let
$$
H=gp\langle G,t|\ A_it=tB_i,\ i\in I\rangle.
$$
Then $H$ is called a group with the stable letter $t$ and the base
group $G$.
\end{definition}

Generally, we may use groups with (many) stable letters $T=\{t\}$.

\noindent{\bf Remark}: Let $H$ be in Definition \ref{d2.1}. P. S.
Novikov (\cite{n52, n54, n55}) called the letter $t$ to be regular
if the subgroup $gp\langle A_i\ | i\in I\rangle,\ gp\langle B_i\ |
i\in I\rangle$ of $G$ are isomorphic by $\varphi: A_i\mapsto B_i,\ \
i\in I$. Thus, Novikov's group ${G}$ with a regular stable letter
$t$ and the base $G$ is exactly an HNN-extension of $G$.

Define the corresponding words relative to a stable letter $t$ by
the above relations:
$$
\mathcal {A}(t)=A_{i_1}^{\pm1}\cdots A_{i_k}^{\pm1},\ \mathcal
{B}(t)=B_{i_1}^{\pm1}\cdots B_{i_k}^{\pm1}.
$$
Moreover, for convenience, we put $\mathcal {A}(t^{-1})=\mathcal
{B}(t)$ and $\mathcal {B}(t^{-1})=\mathcal {A}(t)$. Then, it is
clear that for any $\mathcal {A}(t^{\varepsilon})\in G, \ \mathcal
{A}(t^{\varepsilon})t^{\varepsilon}=t^{\varepsilon}\mathcal
{B}(t^{\varepsilon})$, where $\varepsilon=\pm1$.

Let $G_0\rightarrow G_1\rightarrow \cdots \rightarrow G_n$ be a
tower of groups, where $G_{i+1}$ is a group with some stable letters
and the base group $G_i$ for each $i$. We call such a tower a
Novikov tower. Moreover, if each $G_{i+1}$ is an HNN-extension of
$G_i$, then we call this tower a tower of HNN-extensions or B-tower
(Britton tower, see \cite{br63}).

Let $G_0\rightarrow G_1\rightarrow \cdots \rightarrow G_n$ be a
Novikov tower with $G_0$ free. If $p$ is a stable letter of
$G_{i+1}$, we say the weight of $p$ to be $i+1$. Taking an arbitrary
relation $Ap=pB\ (A,\ B\in G_i)$ from $G_{i+1}$, we can represent it
as follows:
$$
A'xA''p=pB'yB'',
$$
where $x,\ y$ are some stable letters of the highest weight in the
words $A$ and $B$, respectively. We call $x, y$ to be distinguishing
letters of the relation $Ap=pB$. We associate four types of
forbidden subwords in $G_{i+1}$ (see \cite{b94}, \S 6.4):
\begin{eqnarray}\label{Re2.1}
x\mathcal {B}(x)A''p,\ x^{-1}\mathcal {B}(x^{-1})A'^{-1}p,\
y\mathcal {B}(y)B''p^{-1},\ y^{-1}\mathcal {B}(y^{-1})B'^{-1}p^{-1}
\end{eqnarray}
Define the set $C_i$ of words in $G_i\ (0\leq i\leq n)$ as follows
(see \cite{b94, ka, ka2006}).
\begin{enumerate}
\item[(i)]\ $C_0$ is the set of all freely reduced words.
 \item[(ii)]\  Suppose that $C_i$ is defined and the
 algorithm of reducing a word $u\in G_i$ to canonical form $C(u)$ has been defined. For any $w\in C_{i+1},\
 w= w_0p_{1}^{\varepsilon_1}\cdots p_{m}^{\varepsilon_m}w_m$ is freely reduced,
 where $m\geq 0,\ \varepsilon_j=\pm1,\ w_0, w_j\in C_i,\ p_j$ is
stable letter of weight $i+1, \ j=1,2,\cdots,m$,
 and $w$ does not contain the subword as in (\ref{Re2.1}),
  to be more precise, the subwords $xC(\mathcal {B}(x)A'')p$ and so on related to (\ref{Re2.1}).
\end{enumerate}
It is clear that $C_0\subset C_1\subset\cdots\subset C_n$. The
elements of $C_i$ will be called canonical forms.

We describe the algorithm of reducing a word $w=
w_0p_{1}^{\varepsilon_1}\cdots p_{m}^{\varepsilon_m}w_m$ to go to
canonical form by induction on $i$ (see \cite{b94}):
\begin{enumerate}
\item[(a)]\ Reduce the word $w_j$ to canonical form in $G_i$.
\item[(b)]\ Perform all the possible cancelations of letters of weight
$i+1$.
\item[(c)]\ Distinguish one of the occurrences of forbidden subwords
(\ref{Re2.1}) which are related to the stable letter $p=p_i$, where
$i,\ 1\leq i\leq m$, is the minimal index, and eliminate this
subwords by the following rules: \
\begin{eqnarray*}\label{Re2.2}
x\mathcal {B}(x)A''p&=&\mathcal {A}(x)A'^{-1}p B, \ x^{-1}\mathcal
{A}(x)A'^{-1}p=\mathcal {B}(x)A''pB^{-1},\\
 y\mathcal
{B}(y)B''p^{-1}&=&\mathcal {A}(y)B'^{-1}p^{-1}A, \ y^{-1}\mathcal
{A}(y)B'^{-1}p^{-1}=\mathcal {B}(y)B''p^{-1}A^{-1} \ \ \ \ \ \ \ \ \
(2)
\end{eqnarray*}
 where $\mathcal {A}(x),\ \mathcal {B}(x), \ \mathcal {A}(y),\ \mathcal {B}(y)$
are corresponding words (more precisely, by $xC(\mathcal
{B}(x)A'')p=\mathcal {A}(x)A'^{-1}pB$ and so on).
\item[(d)]\ Return to step (a).
\end{enumerate}

\begin{definition}\cite{b94}
The group $G_{n}$ is called a group with standard normal form if for
any $w\in G_{i} \ (0\leq i\leq n),\ w$ can be reducible to canonical
form in a finite steps and the canonical form of $w$ is unique.
\end{definition}

Now we cite the definition of the Gr\"obner-Shirshov basis for the
associative algebra $k\langle X\rangle$ (see \cite{s, b72, b76}).

Let $X$ be a set, $X^*$ the free monoid generated by $X$. We denote
the empty word by $1$, and the length of a word $u$ by $l(u)$. In
general, we set $deg(u)=l(u)$.

A well order $<$ on $X^*$ is monomial if it is compatible with the
multiplication of words, that is, for $u$, $v\in X^*$, we have

$$
u>v \Longrightarrow w_{1}uw_{2}>w_{1}vw_{2},for \ \ all\ \
w_{1},w_{2}\in X^*.
$$

A standard example of monomial order on $X^*$ is the deg-lex order
to compare two words first by degree and then lexicographically,
where $X$ is a linearly ordered set.

Suppose that $X^*$ equipped with a monomial order. Let $f$ and $g$
be two monic polynomials in $k\langle X\rangle$. We denote by
$\bar{f}$ the leading word of $f$. Then, there are two kinds of
compositions:

$(i)$ If \ $w$ is a word such that $w=\bar{f}b=a\bar{g}$ for some
$a,b\in X^*$ with deg$(\bar{f})+$deg$(\bar{g})>$deg$(w)$, then the
polynomial
 $(f,g)_w=fb-ag$ is called the intersection composition of $f$ and
$g$ with respect to $w$.

$(ii)$ If  $w=\bar{f}=a\bar{g}b$ for some $a,b\in X^*$, then the
polynomial $(f,g)_w=f - agb$ is called the inclusion composition of
$f$ and $g$ with respect to $w$.

In the above case, the transformation $f\mapsto(f,g)_w=f-agb$ is
called the elmination of the leading word (ELW) of $g$ in $f$.

Let $S\subset k\langle X\rangle$ be monic. Then the composition
$(f,g)_w$ is called trivial modulo $(S,w)$ if $(f,g)_w=\sum\alpha_i
a_i s_i b_i$, where each $\alpha_i\in k$, $a_i,b_i\in X^{*}, \
s_i\in S$ and $a_i\bar{s_i}b_i<w$. If this is the case, then we
write
$$
(f,g)_w\equiv0\quad mod(S,w).
$$
In general, for $p,q\in k\langle X\rangle$, we write $p\equiv q\quad
mod (S,w)$ which means that $p-q=\sum\alpha_i a_i s_i b_i $, where
each $\alpha_i\in k, \ a_i,b_i\in X^{*}, \ s_i\in S$ and
$a_i\bar{s_i}b_i<w$.

We call the set $S$ with respect to the monomial order $<$ a
Gr\"{o}bner-Shirshov basis in $k\langle X\rangle$ if any composition
of polynomials in $S$ is trivial modulo $S$.

The following lemma was proved by Shirshov \cite{s} for the free Lie
algebras (with deg-lex ordering) in 1962 (see also Bokut
\cite{b72}). In 1976, Bokut \cite{b76} specialized the approach of
Shirshov to associative algebras (see also Bergman \cite{b}). For
commutative polynomials, this lemma is known as the Buchberger's
Theorem (see \cite{bu65}), published in \cite{bu70}.

\begin{lemma}\label{l2.6}
(\cite{b, b76, BoCh1}, Composition-Diamond Lemma) Let $S\subset
k\langle X\rangle$ be a non-empty set  with each $s\in S$ monic and
$<$ a monomial order on $X^*$. Then the following statements are
equivalent:
\begin{enumerate}
\item[(i)]\ $S$ is a Gr\"obner-Shirshov basis in $k\langle X\rangle$.
 \item[(ii)]\ $
f\in Id(S)\Rightarrow\overline{f}=a\overline{s}b$ for some $a,\ b\in
X^*, \ s\in S$.
\item[(iii)]\
$Irr(S)=\{w \in X^*|w\neq a\overline{s}b, \ a,\ b\in X^*, \ s\in
S\}$ is a $k$-linear basis for the factor algebra $k\langle
X|S\rangle$.
\end{enumerate}
\end{lemma}

We now give the definition of the standard Gr\"obner-Shirshov basis
which is associated with the standard normal form on a group (see
\cite{b94}).

\begin{definition}\label{d2.2}
Let $X=Y\dot{\cup} Z$, words $Y^*$ and the letters $Z$ be well
ordered. Suppose that the order on $Y^*$ is monomial. Then, any word
in $X$ has the form $u=u_0z_1\cdots z_ku_k$, where $k\geq 0,\ z_i\in
Z,\ u_i\in Y^*$. Define the weight of the word $u\in X^*$ by
$$
wt(u)=(k, z_1,\cdots, z_k, u_0,\cdots, u_k).
$$
We order $X^*$ as follows.
$$
u>v\Leftrightarrow wt(u)>wt(v).
$$
Then we call the above order the first tower order. Clearly, this
order is a monomial order on $X^*$.
\end{definition}

\begin{definition}\label{d2.3}
Let $X=Y\dot{\cup} \{t^{\pm1}\}$, words $Y^*$ and $\{t^{\pm1}\}$ be
well ordered. Suppose that the order on $Y^*$ is monomial. Then, any
word in $X$ has the form $u=u_0t_1^{\varepsilon_1}\cdots
t_k^{\varepsilon_k}u_k$, where $k\geq 0,\ u_i\in Y^*,\
t_i^{\varepsilon_i}\in \{t^{\pm1}\},\  \varepsilon_i=\pm1$. Define
the weight of the word $u\in X^*$ by
$$
wt(u)=(k_1,k_2, t_1^{\varepsilon_1},\cdots, t_k^{\varepsilon_k},
u_0,\cdots, u_k),
$$
where $k_1$ and  $k_2$  are the numbers of the occurrence of
$t_i^{-1}$ and $t_i$ in $u$, respectively. Now we order $X^*$ as
follows.
$$
u>v\Leftrightarrow wt(u)>wt(v).
$$
Then we call the above order the second tower order. Clearly, this
order is also a monomial order on $X^*$.
\end{definition}

In case $Y=T\dot{\cup} U$ and $Y^*$ are endowed with  one of the
above tower orders, we call the order of words in $X$ the tower
order of words relative to the presentation
$$
X=(T\dot{\cup} U)\dot{\cup} Z.
$$
In general, we can define the tower order of $X$-words relative to
the presentation
$$
X=(\cdots
(X^{(n)}\dot{\cup}X^{(n-1)})\dot{\cup}\cdots)\dot{\cup}X^{(0)},
$$
where $X^{(n)}$-words are endowed by a monomial order.

\begin{definition} Let $G_0\rightarrow G_1\rightarrow \cdots
\rightarrow G_n$ be a Novikov tower as above. If relations
(\ref{Re2.2}) together with trivial relations constitute a
Gr\"{o}bner-Shirshov basis for the group $G_n$ relative to the tower
order, then we call $G_n$ to be a group with the standard
Gr\"{o}bner-Shirshov basis.
\end{definition}

\section{Main result}
In this section, we will deal with two-generator one-relator groups
with depth $k=1,2,3$, respectively and give some examples. Using
Magnus' method, we show that any two-generator one-relator group
with the depth $\leq$ 3 is effectively embedded into a Novikov
tower. Moreover, each group of this tower has effective standard
Gr\"{o}bner-Shirshov basis. Then, from Composition-Diamond lemma, it
follows that each group of this tower has the effective standard
normal form. Also, it follows that this Novikov tower is in fact a
tower of HNN-extensions.

The following is the main result in this paper.
\begin{theorem}\label{t3.1}
Any two-generator one-relator group with the depth $\leq$ 3 is
effectively Magnus embeddable into a tower of HNN-extensions in
which each group has  the effective standard normal form.
\end{theorem}

We prove Theorem \ref{t3.1} step by step.

\subsection{$k=1$}
Let
$$
G=gp\langle x,\ y|x^{n_1}y^{m_1}=1\rangle
$$
be a one-relator group with depth 1. We can suppose that $n_1>0$.
Let
$$
C=gp\langle a, b|(ab^{-m_1})^{n_1}b^{n_1 m_1}=1\rangle
$$
and define a map
$$
\Psi: G\rightarrow C, \ \ x\mapsto ab^{-m_1}, \ y\mapsto b^{n_1}.
$$
Then, $\Psi$ can be extended as a group homomorphism and hence, $G$
can be embedded into $C$. For $i$ an integer, we put $a_i=b^i a
b^{-i}$ and rewrite the defining relation as $a_0 a_{-m_1}\cdots
a_{-(n_1-1)m_1}=1$. Let $A=\{0, -m_1,\cdots, -(n_1-1)m_1 \}, \ \
l=\min{A}$ and $k=\max{A}$.

If $m_1>0$, then $l=-(n_1-1)m_1$ and $k=0$. Let
\begin{eqnarray*}
G_1&=&gp\langle a_i \ (l< i\leq0)| \varnothing\rangle,\\
G_2&=&gp\langle G_1,\  b|\ a_{i+1}b=ba_i \ (l+1\leq i\leq -1),\
a_{l+1}b=b(a_0a_{-m_1}\cdots a_{-(n_1-2)m_1})^{-1}\rangle.
\end{eqnarray*}
Then, $G_2$ is a group with stable letter $b$ and the base $G_1$,
and $C\cong G_2$ by $a\mapsto a_0,\ b\mapsto b$. This means that
$G_1\leq G_2$ is a Novikov tower. Distinguishing letters of defining
relations of $G_2$ can be uniquely defined in all cases except for
the last relation. In the last relation, we define $a_0^{-1},\
a_{l+1}$ to be distinguishing letters.

Now we can get the forbidden subwords of $G_2$ as follows:
$$
a_i^{\varepsilon}b^{\delta} \ (l+1\leq i\leq 0),  \
a_i^{\varepsilon}a_i^{-\varepsilon}, \
b^{\varepsilon}b^{-\varepsilon},
$$
where $\varepsilon, \ \delta=\pm 1$.

Also, we can get the standard rules which would be used to obtain a
Gr\"{o}bner-Shirshov basis of $G_2$ (see the following relations
3.1-3.2).

Let $X=\{a_i^{\pm1}| \ (l< i\leq0)\}\dot{\cup}\{b, b^{-1}\} ,\
a_{l+1}>a_{l+1}^{-1}>a_{l+2}>\cdots>a_0>a_0^{-1}$ and $b^{-1}>b$.
Then we define the first tower order on $X^*$ as Definition
\ref{d2.2}.

At the end, in $G_2$, we have the following standard relations:
\begin{enumerate}
\item[(3.1)]\ $a_i^{\varepsilon}b^{-1}=b^{-1}a_{i+1}^{\varepsilon},\
a_{i+1}^{\varepsilon}b=ba_i^{\varepsilon}$.
\item[(3.2)]\ $a_0^{\varepsilon}b^{-1}=b^{-1}(a_{-m_1+1}\cdots a_{l+1})^{-\varepsilon},\
 a_{l+1}^{\varepsilon}b=b(a_0a_{-m_1}\cdots a_{-(n_1-2)m_1})^{-\varepsilon}$.
\end{enumerate}
Let $S$ consist of the above relations and the trivial relations in
$G_2$. It is easy to check that all compositions in $S$ are trivial.
Thus, with the tower order as above, $S$ is an effective standard
Gr\"{o}bner-Shirshov basis. By Lemma \ref{l2.6},
\begin{eqnarray*}
Irr(S)=&\{&b^na_{i_1}\cdots a_{i_m}|\ n\in\mathbb{Z}, \ m\geq0,\\
&&a_{i_j}\in \{a_i^{\pm 1} \ (l+1\leq i\leq 0)\}, \
a_{i_j}a_{i_{j+1}}\neq a_i^{\varepsilon}a_i^{-\varepsilon}, \
\mbox{and}\  \varepsilon=\pm 1 \ \}
\end{eqnarray*}
is an effective $k$-basis of the algebra $kG_2=k\langle X|S
\rangle$. Since the canonical forms of $G_2$ is $Irr(S)$, $G_2$ is a
group with the effective standard Gr\"{o}bner-Shirshov basis and the
effective standard normal form.

If $m_1<0$, then $k=-(n_1-1)m_1$ and $l=0$. We can use the same
method to get the effective standard normal form of $G_2$.
\begin{example}
Let $G=gp\langle x,\ y|x^2y^2=1\rangle$ and $ C=gp\langle a,\ b\ |\
(ab^{-1})^2b^2=1\rangle. $ Clearly, $C\cong G$ by $x\mapsto
ab^{-1},\ y\mapsto b$. Let $a_i=b^i a b^{-i},\ i=1,2,\cdots,\
H=gp\langle a_0|\varnothing\rangle$ and $G_1=gp\langle H,\ b\
|a_0b=ba_0^{-1}\rangle$. Then $H\leq G_1$ is a Novikov tower and
$C\cong G_1$ by $a\mapsto a_0,\ b\mapsto b$. The forbidden subwords
for group $G_1$ are as follows:
$$
a_0^{\varepsilon}b^{\delta},\ a_0^{\varepsilon}a_0^{-\varepsilon},\
b^{\delta}b^{-\delta}\ \ \ (\varepsilon, \ \delta=\pm 1).
$$
Thus, the effective standard normal words of $G_1$ is $\{b^n a_0^m
|\ n, \ m\in \mathbb{Z}\}$ and $x=a_0 b^{-1}, \ y=b$. Hence, the
effective standard normal words of $G$ is $\{y^n (xy)^m, \
y^n(x^{-1}y^{-1})^{m}x^{-1}|\ n\in\mathbb{Z}, \ m\geq0\}$.
\end{example}

\begin{example}
Let $G=gp\langle x,\ y|x^2y^3=1\rangle$ and $ C=gp\langle a,\ b\ |\
(ab^{-3})^2b^6=1\rangle. $ Clearly, $G\hookrightarrow C$ by
$x\mapsto ab^{-3},\ y\mapsto b^2$ is a monomorphism. Let $a_i=b^i a
b^{-i},\ i=1,2,\cdots,\ H=gp\langle a_i, \ (-2\leq i\leq
0)|\varnothing\rangle$ and $G_1=gp\langle H,\ b\ |a_0b=ba_{-1}, \
a_{-1}b=ba_{-2}, \ a_{-2}b=ba_0^{-1}\rangle $. Then $H\leq G_1$ is a
Novikov tower and $C\cong G_1$ by $a\mapsto a_0,\ b\mapsto b$. The
forbidden subwords for group $G_1$ are as follows:
$$
a_i^{\varepsilon}b^{\delta},\ a_i^{\varepsilon}a_i^{-\varepsilon},\
b^{\delta}b^{-\delta}\ (-2\leq i\leq 0, \ \varepsilon, \ \delta=\pm
1).
$$
Thus, the effective standard normal words of $G_1$ is $\{b^n
a_{i_1}^{\varepsilon_1}\cdots a_{i_m}^{\varepsilon_m}|\ n\in
\mathbb{Z}, \ -2\leq i_j\leq 0,\ \varepsilon_j=\pm 1, \ \mbox{and}\
a_{i_j}^{\varepsilon_j}a_{i_{j+1}}^{\varepsilon_{j+1}}\neq
a_i^{\varepsilon}a_i^{-\varepsilon}\}$.
\end{example}

\subsection{$k=2$}
Let
$$
G=gp\langle x,\ y|x^{n_1}y^{m_1}x^{n_2}y^{m_2}=1 \rangle
$$
be a one-relator group with depth 2. Suppose that $n_1,\ m_i>0$.
There are two cases to consider:

Case 1. $n_1+n_2=0$.

Let $y_i=x^i y x^{-i}$ and rewrite the defining relation as
$y_{n_1}^{m_1} y_{0}^{m_2}=1$. Then we have
\begin{eqnarray*}
H_0&=&gp\langle y_i\ (0\leq i\leq n_1)\ |\ y_{n_1}^{m_1} y_{0}^{m_2}=1\rangle,\\
H_1&=&gp\langle H_0,\ x\ |\ y_{i+1}x=xy_i\ (0\leq i\leq n_1-1)
\rangle,
\end{eqnarray*}
where $H_1\cong G$ and $H_0\leq H_1$ is a Novikov tower. Since $n_1,
\ m_1, \ m_2\neq 0$, similar to the case of depth 1, there exists a
Novikov tower $G_0\leq G_1$ such that $H_0$ can be effectively
embedded into $G_1$, where
\begin{eqnarray*}
G_0&=&gp\langle y_i\ (0< i< n_1), \ a_j\ (l< j\leq0)\ | \
\varnothing \rangle,\\
G_1&=&gp\langle G_0,b\ | \ a_{j+1}b=ba_j\ (l\leq j<0),\
a_{l+1}b=b(a_0
a_{-m_2}\cdots a_{-(m_1-2)m_2})^{-1}\rangle,\\
G_2&=&gp\langle G_1,x\ | \ y_{i+1}x=xy_i\ (0<i<n_1-1), \
y_1x=xb^{m_1}, \ a_0 b^{-m_2}x=xy_{n_1-1}\rangle
\end{eqnarray*}
and $l=-(m_1-1)m_2$. It is clear that $G_0\leq G_1\leq G_2$ is a
Novikov tower and $H_1$ can be embedded into $G_2$. The
distinguishing letters of defining relations of $G_2$ can uniquely
defined in all cases except for the last two relations. Here, we
define $y_1$ and the first $b$ the distinguishing letters in
$y_1x=xb^{m_1}$, and $y_{n_1-1}$ and the last $b^{-1}$ the
distinguishing letters in $a_0 b^{-m_2}x=xy_{n_1-1}$.
  Clearly, $G$ can be embedded into $G_2$. Now we
prove that $G_2$ is a group with effective standard normal form. For
$G_1$ and $G_2$, we have the forbidden subwords:
\begin{eqnarray*}
&&G_1 : \ a_j^{\varepsilon}b^{\delta}, \
a_j^{\varepsilon}a_j^{-\varepsilon}, \
b^{\varepsilon}b^{-\varepsilon}, \
y_i^{\varepsilon}y_i^{-\varepsilon}\\
&&G_2 : \ y_i^{\varepsilon}x^{\delta},\ \ b^{m_1}V(a_j)x^{-1},\
b^{m_2}V(a_j)a_0^{-1} x,\ \ b^{-1}V(a_{j})x^{\varepsilon},\ \
x^{\varepsilon}x^{-\varepsilon}
\end{eqnarray*}
where $l< j\leq0,\ 1\leq i\leq n_1-1,\ \varepsilon,\ \delta=\pm1 $
and $V(X)$ means any group word which is generated by $X$. Similar
to the case of depth 1, $G_1$ is a group with effective standard
normal form.

Let $X=(\{y_i^{\pm1},\
a_j^{\pm1}\}\dot{\cup}\{b^{\pm1}\})\dot{\cup}\{x^{\pm1}\}$. Then we
define a tower order on $X^*$.  We first define deg-lex order on
$\{y_i^{\pm1},\ a_j^{\pm1}\}^*$ and the second tower order on
$(\{y_i^{\pm1},\ a_j^{\pm1}\}\dot{\cup}\{b^{\pm1}\})^*$, and then,
the first tower order on $X^*$. Clearly, this order is a monomial
order on $X^*$.

In $G_2$, we have the following standard relations:
\begin{enumerate}
\item[(3.3)]\ $y_i^{\varepsilon}x^{-1}=x^{-1}y_{i+1}^{\varepsilon},\
y_{i+1}^{\varepsilon}x=xy_i^{\varepsilon},\ (1\leq i\leq n_1-2)$,
\item[(3.4)]\ $y_1^{\varepsilon}x=x(b^{m_1})^{\varepsilon},\
 y_{n_1-1}^{\varepsilon}x^{-1}=x^{-1}(a_0b^{-m_2})^{\varepsilon}$,
 \item[(3.5)]\ $b^{m_2}V(a_j)a_0^{-1}x=V^{(m_2)}(a_j)xy_{n_1-1}^{-1},\
 b^{-1}V(a_j)x=V^{(-1)}(a_j)b^{m_2-1}a_0^{-1}xy_{n_1-1},$
 \item[(3.6)]\ $b^{m_1}V(a_j)x^{-1}=V^{(m_1)}(a_j)x^{-1}y_1,\
 b^{-1}V(a_j)x^{-1}=V^{(-1)}(a_j)b^{m_1-1}a_0^{-1}x^{-1}y_1^{-1}$,
\end{enumerate}
where $V^{(+1)}$ is the result of shifting in $V$ all indices of all
letters with $a_j^{\varepsilon}\mapsto a_{j+1}^{\varepsilon}\
(l+1\leq j\leq -1),\ a_0^{\varepsilon}\mapsto (a_{-m_2+1}\cdots
a_{l+1})^{-\varepsilon}$, and $V^{(-1)}$ with
$a_{j+1}^{\varepsilon}\mapsto a_{j}^{\varepsilon}\ (l+2\leq j\leq
0),\ a_{l+1}^{\varepsilon}\mapsto (a_0a_{-m_2}\cdots
a_{-(m_1-2)m_2})^{-\varepsilon}$. It is clear that $u\in C(V(a_j))$
(the canonical word of $V(a_j)$) if and only if $u$ is a freely
reduced word on $a_j$.

Let $S$ consist of relations $(3.1)-(3.6)$ and the trivial relations
in $G_2$, where we substitute the index set $\{-(n_1-1)m_1\leq i\leq
0\}$ with $\{-(m_1-1)m_2\leq i\leq 0\}$ in the relation
$(3.1)-(3.2)$. Clearly, with the tower order as above, $S$ is an
effective standard Gr\"{o}bner-Shirshov basis. By Lemma \ref{l2.6},
$Irr(S)$ is an effective $k$-basis of the algebra $kG_2=k\langle X|S
\rangle$. Since the canonical forms of $G_2$ is $Irr(S)$, $G_2$ is a
group with the effective standard Gr\"{o}bner-Shirshov basis and the
effective standard normal form.

Case 2. $0\neq n_1+n_2=\alpha, \ 0\neq m_1+m_2=\beta$.

Let
$$
C=gp\langle a, b|(ab^{-\beta})^{n_1}b^{\alpha
m_1}(ab^{-\beta})^{n_2}b^{\alpha m_2}=1\rangle.
$$
Then $G$ can be embedded into $C$ by $x\mapsto ab^{-\beta}, \
y\mapsto b^{\alpha}$. In $C$, the exponent sum of $b$ occurred in
the defining relation is $0$. Let $a_i=b^i a b^{-i}$.

If $n_2>0$, we rewrite the defining relation as
$$
r=a_0 a_{-\beta}\cdots a_{-(n_1-1)\beta}a_{-n_1
\beta+m_1\alpha}\cdots a_{-(n_1+n_2-1)\beta+m_1\alpha}=1.
$$
Let $A =\{ 0, -\beta,\cdots, -(n_1-1)\beta \}, \  B=\{
-n_1\beta+m_1\alpha, \cdots,-(n_1+n_2-1)\beta+\alpha m_1\}, \
l=\min\{A, \ B\}$ and $k=\max\{A, \ B\}$. Then
$$
G_1=gp\langle b, a_i \ (l\leq i\leq k)\ |\ r=1,\ a_{i+1}b=ba_i \
(l\leq i\leq k-1)\rangle
$$
and $C\cong G_1$. Let $G_2=gp\langle a_i\ (l\leq i\leq k)|\
r=1\rangle$. Then $G_2\leq G_1$ is a Novikov tower.

If $A\neq B$, there are seven cases to consider:
\begin{eqnarray*}
 r&=&a_0 a_{-\beta}\cdots (a_{-i\beta}\cdots a_{-(n_1-1)\beta})^2
\cdots a_{-(n_1+n_2-1)\beta+m_1\alpha},\\
 &&a_0 a_{-\beta}\cdots
a_{-i\beta}\cdots a_{-j\beta}\cdots
a_{-(n_1-1)\beta}a_{-i\beta}\cdots
a_{-j\beta},\\
&& a_0 a_{-\beta}\cdots (a_{-i\beta}\cdots a_{-(n_1-1)\beta})^2,\\
&&a_0 a_{-\beta}\cdots a_{-i\beta}\cdots a_{-(n_1-1)\beta}a_{-n_1
\beta+m_1\alpha}\cdots a_o\cdots a_{-i\beta},\\
&& a_0 a_{-\beta}\cdots a_{-(n_1-1)\beta}a_{-n_1 \beta+m_1\alpha}
\cdots a_0 a_{-\beta}\cdots a_{-(n_1-1)\beta}\cdots
a_{-(n_1+n_2-1)\beta+m_1\alpha},\\
&&a_0 a_{-\beta}\cdots a_{-i\beta}\cdots a_{-(n_1-1)\beta}a_0\cdots
a_{-i\beta},\\
&& a_0 a_{-\beta}\cdots a_{-(n_1-1)\beta}a_{-(n_1
\beta)+m_1\alpha}\cdots a_{-(n_1+n_2-1)\beta+m_1\alpha}.
\end{eqnarray*}
If this is the case, it is clear that $G_2$ can be viewed as a free
group. Then, we can use the same method similar to the depth 1 to
get the result.

If $A=B$, then $(a_0 a_{-\beta}\cdots a_{-(n_1-1)\beta})^2=1$. Let
\begin{eqnarray*}
C_0&=&gp\langle c,\ d_i\ (1\leq i\leq n_1-1), \
a_j\ (j\notin A)|\ c^2=1 \rangle,\\
G_0&=&gp\langle C, b|\ a_{i+1}b=ba_i\ (i,i+1\notin A),\
cd_{n_1-1}^{-1}\cdots d_1^{-1}b=ba_{-1},\\
&& \ \ \ \ \ \ \ \ \ \ \ d_ib=ba_{-i\beta-1}\ (i\neq n_1-1),\
a_{-i\beta+1}b=bd_i \rangle.
\end{eqnarray*}
Then $C_0\leq G_0$ is a Novikov tower and  $G_1$ can be embedded
into $G_0$ by $a_0\mapsto cd_{n_1-1}^{-1}\cdots d_1^{-1},\
a_{-i\beta}\mapsto d_i,\ b\mapsto b$. Here, we define $c$ and
$a_{-1}$ the distinguishing letters in $cd_{n_1-1}^{-1}\cdots
d_1^{-1}b=ba_{-1}$. Hence, $G$ can be embedded into $G_0$. The
forbidden subwords for $G_0$ are as follows:
\begin{eqnarray*}
&&a_i^{\pm1}b^{-1},\ a_{i+1}^{\pm1}b\ (i,i+1\notin A),\
a_{-1}^{\pm1}b,\ d_{n_1-1}^{\pm1}b^{-1},\\
&&a_{-i\beta-1}^{\pm1}b^{-1},\ d_i^{\pm1}b^{\pm1}\ (i\neq n_1-1),\
a_{-i\beta+1}^{\pm1}b,\ cd_{n_1-1}^{-1}b,\ cb^{-1}.
\end{eqnarray*}
Let $X=\{c,\ d_i^{\pm1}, \ a_j^{\pm1}|\ 1\leq i\leq n_1-1,\ j\notin
A\}\dot{\cup}\{b^{\pm1}\}$ and
$c>d_{n_1-1}^{-1}>d_{n_1-1}>a_{-1}>\cdots$. Then define the first
tower order on $X^*$. In $G_0$, we have the following standard
relations:
\begin{enumerate}
\item[(3.7)]\ $a_i^{\varepsilon}b^{-1}=b^{-1}a_{i+1}^{\varepsilon},\
a_{i+1}^{\varepsilon}b=ba_i^{\varepsilon},\ (i,\ i+1\notin A),\
a_{-1}^{\varepsilon}b^{-1}=b^{-1}(cd_{n_1-1}^{-1}\cdots
d_1^{-1})^{\varepsilon}$,
\item[(3.8)]\ $c^2=1,\
 cd_{n_1-1}^{-1}b=b(a_{-1}\cdots a_{-(n_1-2)\beta-1}),\
 cb=d_{n_1-1}^{-1}b(a_{-1}\cdots a_{-(n_1-2)\beta-1})^{-1}$,
 \item[(3.9)]\ $a_{-i\beta-1}^{\varepsilon}b^{-1}=b^{-1}d_i^{\varepsilon},\
 d_i^{\varepsilon}b=ba_{-i\beta-1}^{\varepsilon}\ (i\neq n_1-1),\
 a_{-i\beta+1}^{\varepsilon}b=bd_i^{\varepsilon},\
 d_i^{\varepsilon}b^{-1}=b^{-1}a_{-i\beta+1}^{\varepsilon}$.
\end{enumerate}
Let $S$ consist of relations $(3.7)-(3.9)$ and the trivial relations
in $G_0$. Clearly, with the tower order, $S$ is an effective
standard Gr\"{o}bner-Shirshov basis. By Lemma \ref{l2.6}, $Irr(S)$
is an effective $k$-basis of the algebra $kG_0=k\langle X|S
\rangle$. Since the canonical forms of $G_0$ is $Irr(S)$, $G_0$ is a
group with the effective standard Gr\"{o}bner-Shirshov basis and the
effective standard normal form.

If $n_2<0$, we use the same method as above.
 \\

\noindent{\bf Remark}: From the above proof, we know that if there
exists an $a_i$ such that $a_i$ occurs in the relation $r$ only
once, then we can get a Novikov tower of groups with the first group
free.

\subsection{$k=3$}
Let
$$
G=gp\langle x,\
y|x^{n_1}y^{m_1}x^{n_2}y^{m_2}x^{n_3}y^{m_3}=1\rangle
$$
be a one-relator group with depth 3. Suppose that $n_1>0$. There are
two cases to consider:

Case 1. $n_1+n_2+n_3=0$.

Let $y_i=x^i y x^{-i}$ and rewrite the defining relation as
$y_{n_1}^{m_1} y_{n_1+n_2}^{m_2}y_{0}^{m_3}=1$. Let $D=\{0,\ n_1, \
n_1+n_2 \}$ and suppose that $l=\min{D}=0, \ k=\max{D}=n_1+n_2$. By
Magnus' method, we can get a Novikov tower of groups:
\begin{eqnarray*}
H_1&=&gp\langle y_i \ (l\leq i\leq k-1),\ a_j\ (-(m_2-1)m_3< j\leq
0)|\varnothing\rangle\\
&\cong& gp\langle y_i \ (l\leq i\leq k),\ a_j\ (-(m_2-1)m_3\leq
j\leq
0)| \ y_{n_1}^{m_1}a_0\cdots a_{-(m_2-1)m_3}=1\rangle ,\\
H_2&=&gp\langle H_1, \ b| \ a_{j+1}b=ba_j\ (j\neq0),\
a_{-(m_2-1)m_3+1}b=b(y_{n_1}^{m_1}a_0\cdots a_{-(m_2-2)m_3})^{-1}\rangle,\\
H_3&=&gp\langle H_2,\  x|\ y_{i+1}x=xy_i \ (i\neq k-1),\
 a_0b^{-m_3}x=xy_{k-1}, \ y_1x=xb^{m_2}\rangle,
\end{eqnarray*}
where $G$ can be embedded into $H_3$ by $x\mapsto x,\ y\mapsto y_0$.
Here, in $H_2$, we define $a_{-(m_2-1)m_3+1}$ and $a_{0}^{-1}$ the
distinguishing letters in
$a_{-(m_2-1)m_3+1}b=b(y_{n_1}^{m_1}a_0\cdots a_{-(m_2-2)m_3})^{-1}$;
in $H_3$, $y_{k-1}$ and the last $b^{-1}$ the distinguishing letters
in $a_0b^{-m_3}x=xy_{k-1}$, and $y_1$ and the first $b$ the
distinguishing letters in $y_1x=xb^{m_2}$. The forbidden subwords
for the above groups $H_2, H_3$ are as follows:
\begin{eqnarray*}
&&H_2 : \ a_j^{\pm1}b^{-1}, \ \ a_{j+1}^{\pm1}b\ (j\neq0),\
a_{-(m_2-1)m_3+1}^{\pm1}b,\
a_0^{-1}y_{n_1}^{-m_1}b^{-1}, \ a_0b^{-1},\\
&&H_3 : \ y_i^{\pm1}x^{-1},\ y_{i+1}^{\pm1}x\ (1\leq i\leq k-2),\
y_{k-1}^{\pm1}x^{-1},\ y_1^{\pm1}x,\\
&& \ \ \ \ \ \  b^{-1}V(a_{j})x^{\pm1}\ (-(m_2-1)m_3< j\leq 0),\
bV_1(a_j,\ y_{n_1}^{m_1}a_0)\cdots bV_{m_3}(a_j,\
y_{n_1}^{m_1}a_0)a_0^{-1}x,\\
&& \ \ \ \ \ \ bV_1(a_j,\ y_{n_1}^{m_1}a_0)\cdots bV_{m_2}(a_j,\
y_{n_1}^{m_1}a_0)x^{-1}\ (j\neq 0)
\end{eqnarray*}
where for any $1\leq s\leq m_2,m_3,\ V_s^{(s)}$ exists (see
following). Similar to the depth 1, $H_2$ is a group with the
effective standard normal form. Let $X=((\{y_i^{\pm1}|\ l\leq i\leq
k-1\}\dot{\cup}\{a_j^{\pm1}|\ -(m_2-1)m_3< j\leq 0 \})\dot{\cup}\{
b^{\pm1}\})\dot{\cup}\{x^{\pm1}\}$. Define the first tower order on
$(\{y_i^{\pm1}|\ l\leq i\leq k-1\}\dot{\cup}\{a_j^{\pm1}|\
-(m_2-1)m_3< j\leq 0 \})^*$, the second tower order on
$((\{y_i^{\pm1}|\ l\leq i\leq k-1\}\dot{\cup}\{a_j^{\pm1}|\
-(m_2-1)m_3< j\leq 0 \})\dot{\cup}\{ b^{\pm1}\})^*$ and then the
first tower order on $X^*$. In $H_3$, we have the following standard
relations:
\begin{enumerate}
\item[(3.10)]\ $y_i^{\varepsilon}x^{-1}=x^{-1}y_{i+1}^{\varepsilon},\
y_{i+1}^{\varepsilon}x=xy_i^{\varepsilon},\ (1\leq i\leq k-2)$,
\item[(3.11)]\ $y_1^{\varepsilon}x=x(b^{m_2})^{\varepsilon},\
 y_{k-1}^{\varepsilon}x^{-1}=x^{-1}(a_0b^{-m_3})^{\varepsilon},\
a_0^{-1}y_{n_1}^{-m_1}b^{-1}=b^{-1}(a_{-m_3+1}\cdots
a_{-(m_2-1)m_3+1})$,
\item[(3.12)]\ $b^{-1}V(a_{j+1})x=V^{(-1)}(a_{j+1})b^{m_3-1}a_0^{-1}xy_{k-1}$,
\item[(3.13)]\ $b^{-1}V(a_{j+1})x^{-1}=V^{(-1)}(a_{j+1})b^{m_2-1}a_0^{-1}x^{-1}y_1^{-1}$,
\item[(3.14)]\ $bV_1(a_j,\ y_{n_1}^{m_1}a_0)\cdots bV_{m_3}(a_j,\ y_{n_1}^{m_1}a_0)a_0^{-1}x
 =V_1^{(1)}(a_j,\  y_{n_1}^{m_1}a_0)\cdots V_{m_3}^{(m_3)}(a_j,\ y_{n_1}^{m_1}a_0)xy_{k-1}^{-1}$,
\item[(3.15)]\ $bV_1(a_j,\ y_{n_1}^{m_1}a_0)\cdots bV_{m_2}(a_j,\ y_{n_1}^{m_1}a_0)x^{-1}
 =V_1^{(1)}(a_j,\  y_{n_1}^{m_1}a_0)\cdots V_{m_2}^{(m_2)}(a_j,\  y_{n_1}^{m_1}a_0)x^{-1}y_1$,
\item[(3.16)]\ $a_j^{\varepsilon}b^{-1}=b^{-1}a_{j+1}^{\varepsilon},\
a_{j+1}^{\varepsilon}b=ba_j^{\varepsilon},\ (j\neq0),\
a_0b^{-1}=y_{n_1}^{-m_1}b^{-1}(a_{-m_3+1}\cdots
a_{-(m_2-1)m_3+1})^{-1}$,
\item[(3.17)]\ $a_{-(m_2-1)m_3+1}^{\varepsilon}b=b(y_{n_1}^{m_1}a_0\cdots
a_{-(m_2-2)m_3-1})^{-\varepsilon}$,
\end{enumerate}
where the relations $(3.14)-(3.15)$  hold only if $V_s^{(s)}$
exists. Here, $V^{(+1)}$ is the result of shifting in $V$ all
indices of all letters with $a_j^{\varepsilon}\mapsto
a_{j+1}^{\varepsilon}\ (j\neq0),\
y_{n_1}^{m_1}a_0^{\varepsilon}\mapsto (a_{-m_3+1}\cdots
a_{-(m_2-1)m_3+1})^{-\varepsilon}$ and $V^{(-1)}$ with
$a_{j+1}^{\varepsilon}\mapsto a_{j}^{\varepsilon}\ (j\neq
-(m_2-1)m_3+1),\ a_{-(m_2-1)m_3+1}^{\varepsilon}\mapsto
(y_{n_1}^{m-1}a_0a_{-m_2}\cdots a_{-(m_2-2)m_3})^{-\varepsilon}$.
Since $\{a_{j+1}\}\ \ (\{a_j,\ y_{n_1}^{m_1}a_0\})$ freely generate
a subgroup of $H_1$, $C(V(a_{j+1}))\ \ (C(V(a_j,\
y_{n_1}^{m_1}a_0))$ is the freely reduced word on $\{a_{j+1}\}\ \
(\{a_j,\ y_{n_1}^{m_1}a_0\})$.

Let $S$ consist of relations $(3.10)-(3.17)$ and the trivial
relations in $H_3$. Clearly, with the tower order mentioned as
above, $S$ is an effective standard Gr\"{o}bner-Shirshov basis. By
Lemma \ref{l2.6}, $Irr(S)$ is an effective $k$-basis of the algebra
$kH_3=k\langle X|S \rangle$. Since the canonical forms of $H_3$ is
$Irr(S)$, $H_3$ is a group with the effective standard
Gr\"{o}bner-Shirshov basis and the effective standard normal form.

Case 2. $0\neq n_1+n_2+n_3=\alpha, \ 0\neq m_1+m_2+m_3=\beta$.

We may assume that all the powers are positive. For other cases, we
use a similar way to prove the result. Let
\begin{eqnarray*}
G_1&=&gp\langle a,\ b |\ (ab^{-\beta})^{n_1}b^{\alpha
m_1}(ab^{-\beta})^{n_2}b^{\alpha m_2}(ab^{-\beta})^{n_3}b^{\alpha
m_3}=1\rangle,\\
G_2&=&gp\langle b, a_i\ (l\leq i\leq k)|\ r=1,\ a_{i+1}b=ba_i\
(l\leq i<k)\rangle,
\end{eqnarray*}
where $l=min\{A,\ B,\ C\},\ k=max\{A,\ B,\ C\}$,
\begin{eqnarray*}
r&=&a_0a_{-\beta}\cdots
a_{-(n_1-1)\beta}a_{-n_1\beta+m_1\alpha}\cdots
a_{-(n_1+n_2-1)\beta+m_1\alpha}\cdot\\
&& \ \ a_{-(n_1+n_2)\beta+(m_1+m_2)\alpha}\cdots
a_{-(n_1+n_2+n_3-1)\beta+(m_1+m_2)\alpha},\\
A&=&\{0, -\beta,\cdots,-(n_1-1)\beta\},\\
B&=&\{-n_1\beta+m_1\alpha,\cdots,-(n_1+n_2-1)\beta+m_1\alpha \},\\
C&=&\{-(n_1+n_2)\beta+(m_1+m_2)\alpha,\cdots,-(n_1+n_2+n_3-1)\beta+(m_1+m_2)\alpha\}.
\end{eqnarray*}
Then $G\hookrightarrow G_1\simeq G_2$. Let $H_1=gp\langle a_i\
(l\leq i\leq k)|\ r=1\rangle, \ H_2=gp\langle H_1,\ b|\
a_{i+1}b=ba_i\ (l\leq i<k)\rangle$. It is clear that $H_1\leq H_2$
is a Novikov tower and $G_2\cong H_2$.

If there exists an $i$ such that $i$ is only in one of the sets $
A,\ B,\ \mbox{or}\ C$, then we get the result by the Remark in 3.2.

If $A=B\cup C$, then $r$ has four possibilities:
\begin{eqnarray*}
r&=&a_0a_{-\beta}\cdots a_{-(n_1-1)\beta}a_{-s\beta}\cdots
a_{-(n_1-1)\beta}a_0\cdots a_{-(s-1)\beta},\\
&&(a_0a_{-\beta}\cdots a_{-(n_1-1)\beta})^2,\\
&&a_0a_{-\beta}\cdots a_{-(n_1-1)\beta}a_{-s\beta}\cdots
a_{-(n_1-1)\beta}a_0\cdots a_{-t\beta}\ (0\leq s\leq t\leq n_1-1),\\
&&a_0a_{-\beta}\cdots a_{-(n_1-1)\beta}a_o\cdots
a_{-t\beta}a_{-s\beta}\cdots a_{-(n_1-1)\beta}\ (0\leq s\leq t\leq
n_1-1).
\end{eqnarray*}
We only consider the first case. Other cases can be similarly
proved.

Similar to the depth 2, we can get a Novikov tower of groups:
\begin{eqnarray*}
C_0&=&gp\langle c,\ a_i\ (i\notin A),\ d_j(1\leq j\leq
n_1-2)|\varnothing \rangle,\\
C_1&=&gp\langle C_0,\ d|cd=dc^{-1}
\rangle,\\
C_2&=&gp\langle C_1,\  b|\ a_{i+1}b=ba_i\ (i,\ i+1\notin A),\
cd^{-1}d_{s-1}^{-1}\cdots d_1^{-1}b=ba_{-1},\\
&&a_{-i\beta+1}b=bd_i,\ d_ib=ba_{-i\beta+1}\ (1\leq i\leq
s-1),\\
&&a_{-j\beta+1}b=bd_{j-1}\ (s+1\leq j\leq n_1-1),\
d_{j-1}b=ba_{-j\beta-1}\ (s+1\leq
j\leq n_1-2),\\
&&cd_{n_1-2}^{-1}\cdots d_s^{-1}b=ba_{-1}\cdots a_{-s\beta-1},\
a_{-s\beta+1}b=bdd_{n_1-2}^{-1}\cdots d_s^{-1} \rangle,
\end{eqnarray*}
where $G_2\hookrightarrow C_2$ by $a_i\mapsto a_i\ (i\notin A),\
a_0\mapsto cd^{-1}d_{s-1}^{-1}\cdots d_1^{-1},\ a_{-i\beta}\mapsto
d_i\ (1\leq i\leq s-1),\ a_{-j\beta}\mapsto d_{j-1}\ (s+1\leq j\leq
n_1-1),\ a_{-s\beta}\mapsto dd_{n_1-2}^{-1}\cdots d_s^{-1}$. Here,
in $C_2$, we define $a_{-1}$ and $d^{-1}$ in
$cd^{-1}d_{s-1}^{-1}\cdots d_1^{-1}b=ba_{-1},\ \ a_{-s\beta-1}$ and
$d_{n_1-2}^{-1}$ in $cd_{n_1-2}^{-1}\cdots d_s^{-1}b=ba_{-1}\cdots
a_{-s\beta-1}$, and $a_{-s\beta+1}$ and $d$ in
$a_{-s\beta+1}b=bdd_{n_1-2}^{-1}\cdots d_s^{-1}$ the distinguishing
letters. The forbidden subwords for $C_1$ and $C_2$ are:
\begin{eqnarray*}
&&C_1 : \ c^{\varepsilon}d^{\delta},\\
&&C_2 : \ a_i^{\varepsilon}b^{\delta},\
d_{i}^{\varepsilon}b^{\delta}\ (i\neq n_1-2),\
d_{n_1-2}^{\varepsilon}b^{-\varepsilon},\ d_{n_1-2}c^{-1}b,\
dc^{m}d_{n_1-2}^{-1}b^{-1},\ d^{-1}c^mb^{\varepsilon},\
dc^{m}c^{-1}b.
\end{eqnarray*}
Let $X=(\{c^{\pm1},\ a_i^{\pm1},\ d_j^{\pm1}|\ i\notin A,\ 1\leq
j\leq n_1\}\dot{\cup} \{d^{\pm1}\})\dot{\cup}\{b^{\pm1}\}$ and
$d_{n_1-2}>d_{n_1-2}^{-1}>\cdots>a_i>a_i^{-1}>c>c^{-1}$. Then define
the first tower order on $X^*$. In $C_2$, we have the following
standard relations:
\begin{enumerate}
\item[(3.18)]\ $c^{\varepsilon}d^{\delta}=d^{\delta}c^{-\varepsilon},\
a_{i+1}^{\varepsilon}b=ba_i^{\varepsilon},\
a_i^{\varepsilon}b^{-1}=b^{-1}a_{i+1}^{\varepsilon} (i,\ i+1\notin
A),\ a_{-1}^{\varepsilon}b^{-1}=b^{-1}(cd^{-1}d_{s-1}^{-1}\cdots
d_1^{-1})^{\varepsilon}$,
\item[(3.19)]\ $d^{-1}c^mb=c^{-m-1}ba_{-1}\cdots a_{-(s-1)\beta-1},\
(dc^m)c^{-1}b=c^{-m}b(a_{-1}\cdots a_{-(s-1)\beta-1})^{-1}$,
\item[(3.20)]\ $d_i^{\varepsilon}b^{-1}=b^{-1}a_{-i\beta+1}^{\varepsilon},\
a_{-i\beta+1}^{\varepsilon}b=bd_i^{\varepsilon},\
d_i^{\varepsilon}b=ba_{-i\beta-1}^{\varepsilon},\
a_{-i\beta-1}^{\varepsilon}b^{-1}=b^{-1}d_i^{\varepsilon}\ (1\leq
i\leq s-1)$,
\item[(3.21)]\ $d_{j-1}^{\varepsilon}b^{-1}=b^{-1}a_{-j\beta+1}^{\varepsilon},\
a_{-j\beta+1}^{\varepsilon}b=bd_{j-1}^{\varepsilon}\ (s+1\leq j\leq
n_1-1)$,
 \item[(3.22)]\ $
d_{j-1}^{\varepsilon}b=ba_{-j\beta-1}^{\varepsilon},\
a_{-j\beta-1}^{\varepsilon}b^{-1}=b^{-1}d_{j-1}^{\varepsilon}\
(s+1\leq j\leq n_1-2)$,
\item[(3.23)]\ $d_{n_1-2}c^{-1}b=b(a_{-1}\cdots
a_{-(n_1-2)\beta-1})^{-1},\
dc^md_{n_1-2}^{-1}b^{-1}=c^{-m}b^{-1}(a_{-s\beta+1}\cdots
a_{-(n_1-2)\beta+1})$,
 \item[(3.24)]\ $d^{-1}c^mb^{-1}=c^{-m}d_{n_1-2}^{-1}b^{-1}(a_{-s\beta+1}\cdots
a_{-(n_1-2)\beta+1})^{-1},\ d_{n_1-2}^{-1}b=c^{-1}b(a_{-1}\cdots
a_{-(n_1-2)\beta-1})$,
 \item[(3.25)]\ $a_{-s\beta-1}^{\varepsilon}b^{-1}=b^{-1}(dd_{n_1-2}^{-1}\cdots
d_s^{-1})^{\varepsilon},\
a_{-s\beta+1}^{\varepsilon}b=b(dd_{n_1-2}^{-1}\cdots
d_s^{-1})^{\varepsilon}$.
\end{enumerate}

Let $S$ consist of relations $(3.18)-(3.25)$ and the trivial
relations in $C_2$. Clearly, with the tower order as above, $S$ is
an effective standard Gr\"{o}bner-Shirshov basis. By Lemma
\ref{l2.6}, $Irr(S)$ is an effective $k$-basis of the algebra
$kC_2=k\langle X|S \rangle$. Since the canonical forms of $C_2$ is
$Irr(S)$, $C_2$ is a group with the effective standard
Gr\"{o}bner-Shirshov basis and the effective standard normal form.

When we replace $A$ with $B$ or $C$, using the same method, we can
get the result.

 \ \

By the above proof, each group in the Novikov tower has effective
standard normal form. From this it follows that each Novikov tower
is exactly a tower of HNN-extensions. Thus, the proof of Theorem
\ref{t3.1} is completed.
 \\

\section{An example}
In this section, we will give an example to show how to deal with
some general cases for one-relator groups.

\begin{example}
 $G=gp\langle a,\ b,\ c,\ t_1|(a^2c^2t_1c^{-4}b^3t_1^{-1}
 c^2t_1^{-2}a^{-3}b^{-2}t_1^{2}ab^2t_1c^{-2}t_1^{-1})^{r}=1\rangle\ (r\geq 2)$.

Let $a^{(0)}_{(i)}=t_1^{i}at_1^{-i},\ b_{(i)}=t_1^{i}bt_1^{-i},\
c_{(i)}=t_1^{i}ct_1^{-i}$ and rewrite the defining relation as
$((a^{(0)}_{(0)})^2c_{(0)}^2c_{(1)}^{-4}b_{(1)}^{3}
c_{(0)}^2(a^{(0)}_{(-2)})^{-3}b_{(-2)}^{-2}a^{(0)}_{(0)}b_{(0)}^2c_{(1)}^{-2})^r=1$.
Then, we get
\begin{eqnarray*}
G_2&=&gp\langle a^{(0)}_{(i)},\ b_{(j)},\ c_{(0)},\ c_{(1)},\
(-2\leq i\leq 0,\ -2\leq j\leq 1) |\\
&&\ \ \ \ \ \ \ \ \
((a^{(0)}_{(0)})^2c_{(0)}^2c_{(1)}^{-4}b_{(1)}^{3}
c_{(0)}^2(a^{(0)}_{(-2)})^{-3}b_{(-2)}^{-2}a^{(0)}_{(0)}b_{(0)}^2c_{(1)}^{-2})^r=1
 \rangle,\\
G_1&=&gp\langle G_2,\ t_1\ |\ a^{(0)}_{(i)}t_1=t_1a^{(0)}_{(i-1)},
b_{(j)}t_1=t_1b_{(j-1)},\ c_{(1)}t_1=t_1c_{(0)},\\
&&\ \ \ \ \ \ \ \ \ \ \ \ \ (-2\leq i-1\leq -1,\ -2\leq j-1\leq 0)
\rangle.
\end{eqnarray*}
Clearly, $G\cong G_1$. Let $a^{(0)}_{(0)}=a^{(01)}_{(0)}t_2,\
a^{(0)}_{(-2)}=t_2$ and rewrite the defining relation of $G_2$ as
$(a^{(01)}_{(0)}t_2a^{(01)}_{(0)}t_2c_{(0)}^2c_{(1)}^{-4}b_{(1)}^{3}
c_{(0)}^2t_2^{-3}b_{(-2)}^{-2}a^{(01)}_{(0)}t_2b_{(0)}^2c_{(1)}^{-2})^r=1$.
Let $a^{(01)}_{(0i)}=t_2^{i}a^{(01)}_{(0)}t_2^{-i},\
b_{(ji)}=t_2^{i}b_{(j)}t_2^{-i},\ c_{(ji)}=t_2^{i}c_{(j)}t_2^{-i},\
a^{(0)}_{(j)}=a^{(00)}_{(j0)}$. Let
\begin{eqnarray*}
G_3&=&gp\langle a^{(01)}_{(j_1i)},\ b_{(j_2i)},\ c_{(0i)},\
c_{(1i)},\ (-2\leq j_1\leq 0,\ -2\leq j_2\leq 1) |\
a^{(01)}_{(01)}t_2=t_2a^{(01)}_{(00)},\\
&&a^{(01)}_{(00)}t_2=t_2a^{(01)}_{(0-1)},\
b_{(-20)}t_2=t_2b_{(-2-1)},\ b_{(12)}t_2=t_2b_{(11)},\
b_{(11)}t_2=t_2b_{(10)},\\
&&c_{(12)}t_2=t_2c_{(11)},\ c_{(11)}t_2=t_2c_{(10)},\
c_{(02)}t_2=t_2c_{(01)},\
c_{(01)}t_2=t_2c_{(00)},\\
&&(a^{(01)}_{(00)}a^{(01)}_{(01)}c_{(02)}^2c_{(12)}^{-4}b_{(12)}^{3}
c_{(02)}^2b_{(-2-1)}^{-2}a^{(01)}_{(0-1)}b_{(00)}^2c_{(10)}^{-2})^r=1\rangle.
\end{eqnarray*}
Let $a^{(01)}_{(00)}=a^{(011)}_{(00)}t_3^{-1},\ a^{(01)}_{(01)}=t_3$
and rewrite the last defining relation of $G_3$ as
$$
(a^{(011)}_{(000)}c_{(020)}^2c_{(120)}^{-4}b_{(120)}^{3}
c_{(020)}^2b_{(-2-10)}^{-2}a^{(01)}_{(0-10)}b_{(000)}^2c_{(100)}^{-2})^r=1.
$$
Then, just like $G_3$, we have $G_4$. Let
$a^{(011)}_{(000)}=a^{(0111)}_{(000)}t_4^{-4},\ c_{(020)}=t_4$ and
rewrite above relation as
$$
(a^{(0111)}_{(000)}t_4^{-2}c_{(120)}^{-4}b_{(120)}^{3}
t_4^2b_{(-2-10)}^{-2}a^{(010)}_{(0-10)}b_{(000)}^2c_{(100)}^{-2})^r=1.
$$
Let $a^{(0111)}_{(000i)}=t_4^{i}a^{(0111)}_{(000)}t_4^{-i},\
b_{(jkli)}=t_4^{i}b_{(jkl)}t_4^{-i},\
c_{(jkli)}=t_4^{i}c_{(jkl)}t_4^{-i}$. Then, we have
\begin{eqnarray*}
&&G_5=gp\langle X,\ t_4\ |\ b_{(1200)}t_4=t_4b_{(120-1)},\
b_{(120-1)}t_4=t_4b_{(120-2)},\ c_{(1200)}t_4=t_4c_{(120-1)},
\\
&& \ \ \ \ \ \ \ c_{(120-1)}t_4=t_4c_{(120-2)},\
(a^{(0111)}_{(0000)}c_{(120-2)}^{-4}b_{(120-2)}^{3}
b_{(-2-100)}^{-2}a^{(0100)}_{(0-100)}b_{(0000)}^2c_{(1000)}^{-2})^r=1\rangle,
\end{eqnarray*}
where $X=\{ a^{(0111)}_{(0000)},\ a^{(0000)}_{(-1000)},\
a^{(0100)}_{(0-100)},\ b_{(120i)},\  b_{(1100)},\ b_{(i000)},\
b_{(-2-100)},\  c_{(jk00)},\ c_{(120i)}| -2\leq i\leq 0,\
j,k=0,1\}$.
 Let
$a^{(0111)}_{(0000)}=a^{(01111)}_{(0000)}t_5^{4},\ c_{(120-2)}=t_5$
and rewrite above relation as
$$
(a^{(01111)}_{(00000)}b_{(120-20)}^{3}
b_{(-2-1000)}^{-2}a^{(01000)}_{(0-1000)}b_{(00000)}^2c_{(10000)}^{-2})^r=1.
$$
Let $a^{(01^4)}_{(0^5)}=a^{(01^5)}_{(0^6)}t_6^{-3},\
b_{(120-20)}=t_6$ and rewrite above relation as
$$
(a^{(01^5)}_{(0^6)}
b_{(-2-10^4)}^{-2}a^{(010^4)}_{(0-10^4)}b_{(0^6)}^2c_{(10^5)}^{-2})^r=1,
$$
where, for example, $(01^4)=(01111)$. Let
$a^{(01^5)}_{(0^6)}=a^{(01^6)}_{(0^7)}t_7^{2},\ b_{(-2-10^4)}=t_7$
and rewrite above relation as
$$
(a^{(01^6)}_{(0^7)}
a^{(010^5)}_{(0-10^5)}b_{(0^7)}^2c_{(10^6)}^{-2})^r=1.
$$
Repeating this process, we can get a Novikov tower:
\begin{eqnarray*}
H_1&=&gp\langle X_1|\ (a^{(01^9)}_{(0^{10})})^r=1
 \rangle,\\
H_2&=&gp\langle H_1,\ t_3,\ t_4,\ t_7\ |\ \
b_{(120^8)}t_4=t_4b_{(120-10^6)},\ b_{(120-10^6)}t_4=t_4t_6,\\
&&c_{(120^8)}t_4=t_4c_{(120-10^6)}, \
c_{(120-10^6)}t_4=t_4t_5 \rangle,\\
H_3&=&gp\langle H_2,\ t_2\ |\ b_{(-20^9)}t_2=t_2t_7,\
b_{(120^8)}t_2=t_2b_{(110^8)},\ b_{(110^8)}t_2=t_2b_{(10^9)},\\
&& c_{(110^8)}t_2=t_2t_{10},\ t_{4}t_2=t_2c_{(010^8)},\
c_{(010^8)}t_2=t_2c_{(0^{10})},
c_{(120^8)}t_2=t_2c_{(110^8)},\\
&&t_3t_2=t_2a^{(01^9)}_{(0^{10})}t_{10}^2t_9^{-2}t_8t_7^{2}t_6^{-3}t_5^4t_4^{-4}t_3^{-1},\
a^{(01^9)}_{(0^{10})}t_{10}^2t_9^{-2}t_8t_7^{2}t_6^{-3}t_5^4t_4^{-4}t_3^{-1}t_2=t_2t_8^{-1}
\rangle,\\
H_4&=&gp\langle H_3,\ t_1\ |\ t_2t_8^{-1}t_1
=t_1a^{(0^{10})}_{(-10^9)},\ a^{(0^{10})}_{(-10^9)}t_1=t_1t_2,\\
&&b_{(10^9)}t_1=t_1t_9,\ t_9t_1=t_1b_{(-10^9)},\
b_{(-10^9)}t_1=t_1b_{(-20^9)},\ t_{10}t_1=t_1c_{(0^{10})}\rangle,
\end{eqnarray*}
where $X_1=\{a^{(01^9)}_{(0^{10})},\ a^{(0^{10})}_{(-10^9)},\
b_{(10^9)},\ b_{(110^8)},\ b_{(120^8)},\ b_{(120-10^6)},\
b_{(-10^9)},\ b_{(-20^9)},\ c_{(0^{10})},\ c_{(010^8)}$,\\ $
c_{(110^8)},\ c_{(120^8)},\ c_{(120-10^6)},\ t_i|5\leq i\leq10,\
i\neq 7\}$.

 The standard relations for groups $H_1,H_2,H_3,H_4$ are as follows:
\begin{enumerate}
\item[(4.1)]\  $(a^{(01^9)}_{(0^{10})})^{(r-i)\varepsilon}=(a^{(01^9)}_{(0^{10})})^{i\varepsilon},\
(0\leq i\leq \lceil r/2\rceil)$,
\item[(4.2)]\ $b_{(120^8)}^{\varepsilon}t_4=t_4b_{(120-10^6)}^{\varepsilon},\
 b_{(120-10^6)}^{\varepsilon}t_4^{-1}=t_4^{-1}b_{(120^8)}^{\varepsilon},\
 b_{(120-10^6)}^{\varepsilon}t_4=t_4t_6^{\varepsilon}$,
\item[(4.3)]\ $t_6^{\varepsilon}t_4^{-1}=t_4^{-1}b_{(120-10^6)}^{\varepsilon},\
c_{(120^8)}^{\varepsilon}t_4=t_4c_{(120-10^6)}^{\varepsilon},\
 c_{(120-10^6)}^{\varepsilon}t_4^{-1}=t_4^{-1}c_{(120^8)}^{\varepsilon}$,
\item[(4.4)]\ $c_{(120-10^6)}^{\varepsilon}t_4=t_4t_5^{\varepsilon},\
t_5^{\varepsilon}t_4^{-1}=t_4^{-1}c_{(120-10^6)}^{\varepsilon}$,
\item[(4.5)]\ $b_{(-20^9)}^{\varepsilon}t_2=t_2t_7^{\varepsilon},\
t_7^{\varepsilon}t_2^{-1}=t_2^{-1}b_{(-20^9)}^{\varepsilon},\
b_{(120^8)}^{\varepsilon}t_2=t_2b_{(110^8)}^{\varepsilon}$,
\item[(4.6)]\ $b_{(110^8)}^{\varepsilon}t_2^{-1}=t_2^{-1}b_{(120^8)}^{\varepsilon},\
b_{(110^8)}^{\varepsilon}t_2=t_2b_{(10^9)}^{\varepsilon},\
b_{(10^9)}^{\varepsilon}t_2^{-1}=t_2^{-1}b_{(110^8)}^{\varepsilon}$,
\item[(4.7)]\ $c_{(120^8)}^{\varepsilon}t_2=t_2c_{(110^8)}^{\varepsilon},\
c_{(110^8)}^{\varepsilon}t_2^{-1}=t_2^{-1}c_{(120^8)}^{\varepsilon},\
c_{(110^8)}^{\varepsilon}t_2=t_2t_{10}^{\varepsilon}$,
\item[(4.8)]\ $t_{10}^{\varepsilon}t_2^{-1}=t_2^{-1}c_{(110^8)}^{\varepsilon},\
c_{(010^8)}^{\varepsilon}t_2=t_2c_{(0^{10})}^{\varepsilon} ,\
c_{(010^8)}^{\varepsilon}t_2^{-1}=t_2^{-1}t_{4}^{\varepsilon}$,
\item[(4.9)]\ $c_{(0^{10})}^{\varepsilon}t_2^{-1}=t_2^{-1}c_{(010^8)}^{\varepsilon},\
t_{4}V(Y)t_2=V'(Y')t_2c_{(010^8)},\
t_{4}^{-1}V'(Y')t_2=V(Y)t_2c_{(010^8)}^{-1}$,
\item[(4.10)]\ $t_3^{\varepsilon}t_2=
t_2(a^{(01^9)}_{(0^{10})}t_{10}^2t_9^{-2}t_8t_7^{2}t_6^{-3}t_5^4t_4^{-4}t_3^{-1})^{\varepsilon}$,
\item[(4.11)]\
$t_3^{-1}t_2^{-1}=(a^{(01^9)}_{(0^{10})}t_{10}^2t_9^{-2}t_8t_7^{2}t_6^{-3}t_5^4t_4^{-4})^{-1}t_2^{-1}t_3$,
\item[(4.12)]\ $(a^{(01^9)}_{(0^{10})}t_{10}^2t_9^{-2}t_8t_7^{2}t_6^{-3}t_5^4t_4^{-4}t_3^{-1})^{-1}t_2^{-1}=
t_2^{-1}t_{3}^{-1}$,
\item[(4.13)]\ $t_8^{\varepsilon}t_2^{-1}=t_2^{-1}(a^{(01^9)}_{(0^{10})}t_{10}^2t_9^{-2}
t_8t_7^{2}t_6^{-3}t_5^4t_4^{-4}t_3^{-1})^{\varepsilon}$,
\item[(4.14)]\ $t_7t_6^{-3}t_5^4t_2=(a^{(01^9)}_{(0^{10})}t_{10}^2
t_9^{-2}t_8t_7)^{-1}t_{2}t_8^{-1}a^{(01^9)}_{(0^{10})}t_{10}^2
t_9^{-2}t_8t_7^{2}t_6^{-3}t_5^4t_4^{-4}t_3^{-1}c_{(010^8)}^{4}$,

\item[(4.15)]\ $(a^{(01^9)}_{(0^{10})}t_{10}^2t_9^{-2}
t_8t_7^{2})^{-1}t_2=
t_6^{-3}t_5^4t_2(a^{(01^9)}_{(0^{10})}t_{10}^2t_9^{-2}t_8t_7^{2}
t_6^{-3}t_5^4t_4^{-4}t_3^{-1}c_{(010^8)}^{4})^{-1}t_8$,

\item[(4.16)]\ $(a^{(0^{10})}_{(-10^9)})^{\varepsilon}t_1=t_1t_2^{\varepsilon},\
t_2V(Y_1)t_1^{-1}=V'(Y'_1)t_1^{-1}(a^{(0^{10})}_{(-10^9)}),\
t_2^{-1}V'(Y'_1)t_1^{-1}=V(Y_1)t_1^{-1}(a^{(0^{10})}_{(-10^9)})^{-1}$,

\item[(4.17)]\ $b_{(10^9)}^{\varepsilon}t_1=t_1t_9^{\varepsilon},\
t_9^{\varepsilon}t_1^{-1}=t_1^{-1}b_{(10^9)}^{\varepsilon},\
t_9^{\varepsilon}t_1=t_1b_{(-10^9)}^{\varepsilon}$,
\item[(4.18)]\ $b_{(-10^9)}^{\varepsilon}t_1^{-1}=t_1^{-1}t_9^{\varepsilon},\
b_{(-10^9)}^{\varepsilon}t_1=t_1b_{(-20^9)}^{\varepsilon},\
b_{(-20^9)}^{\varepsilon}t_1^{-1}=t_1^{-1}b_{(-10^9)}^{\varepsilon}$,
\item[(4.19)]\ $t_{10}^{\varepsilon}t_1=t_1c_{(0^{10})}^{\varepsilon},\
c_{(0^{10})}^{\varepsilon}t_1^{-1}=t_1^{-1}t_{10}^{\varepsilon},\
(a^{(0^{10})}_{(-10^9)})^{\varepsilon}t_1^{-1}=t_1^{-1}(t_2t_8^{-1})^{\varepsilon}$,
\item[(4.20)]\ $t_2V(Y_1)t_8^{-1}t_1=V'(Y'_1)t_1a^{(0^{10})}_{(-10^9)},\
t_2^{-1}V'(Y'_1)t_1
=V(Y_1)t_8^{-1}t_1(a^{(0^{10})}_{(-10^9)})^{-1}$,
\end{enumerate}
where
\begin{eqnarray*}
Y&=&\{b_{(120-10^6)},\ c_{(120-10^6)},\ t_5,\ t_6\},\\
Y'&=&\{b_{(120-10^6)},\ c_{(120-10^6)},\ b_{(120^8)},\
c_{(120^8)}\},\\
Y_1&=&\{t_7,\ t_8,\ t_{10},\ b_{(110^8)},\ b_{(10^9)},\
c_{(110^8)},\ c_{(010^8)},\
c_{(0^{10})},\ a^{(01^9)}_{(0^{10})}t_{10}^2t_9^{-2}t_8t_7^{2}t_6^{-3}t_5^4t_4^{-4}t_3^{-1}\},\\
Y'_1&=&\{t_3,\ t_4,\ b_{(-20^9)},\ b_{(120^8)},\ b_{(110^8)},\
c_{(120^8)},\  c_{(110^8)},\ c_{(010^8)},\
a^{(01^9)}_{(0^{10})}t_{10}^2t_9^{-2}t_8t_7^{2}t_6^{-3}t_5^4\},
\end{eqnarray*}
and $V(Y)\leftrightarrow V'(Y')$ by $b_{(120-10^6)}\leftrightarrow
b_{(120^8)},\ c_{(120-10^6)}\leftrightarrow c_{(120^8)},\
t_5\leftrightarrow c_{(120-10^6)},\ t_6\leftrightarrow
b_{(120-10^6)},\ V(Y_1)\leftrightarrow V'(Y'_1)$ by
$t_7\leftrightarrow b_{(-20^9)},\ t_8\leftrightarrow
(a^{(01^9)}_{(0^{10})}t_{10}^2t_9^{-2}t_8t_7^{2}t_6^{-3}t_5^4t_4^{-4}t_3^{-1})^{-1},\
t_{10}\leftrightarrow c_{(110^8)},\ b_{(110^8)}\leftrightarrow
b_{(120^8)},\ b_{(10^9)}\leftrightarrow b_{(110^8)},\
c_{(110^8)}\leftrightarrow c_{(120^8)},\ c_{(010^8)}\leftrightarrow
t_4,\ c_{(0^{10})}\leftrightarrow c_{(010^8)}$,
$a^{(01^9)}_{(0^{10})}t_{10}^2t_9^{-2}t_8t_7^{2}t_6^{-3}t_5^4\leftrightarrow
t_8^{-1}a^{(01^9)}_{(0^{10})}t_{10}^2t_9^{-2}t_8t_7^{2}t_6^{-3}t_5^4t_4^{-4}t_3^{-1}c_{(010^8)}^{4},\
a^{(01^9)}_{(0^{10})}t_{10}^2t_9^{-2}t_8t_7^{2}t_6^{-3}t_5^4t_4^{-4}t_3^{-1}\leftrightarrow
t_3$. The left part of the relations $(4.2)-(4.4)$ are the forbidden
subwords for $H_2$, $(4.5)-(4.15)$ for $H_3$ and $(4.16)-(4.20)$ for
$H_4$. Since $Y\ (Y')$ freely generates a subgroup of $H_1$,
$C(V(Y))\ (C(V'(Y')))$ is the freely reduced word on $Y\ (Y')$.
$u\in C(V(Y_1))$ if and only if $u$ is reduced on $\{Y_1\backslash
\{a^{(01^9)}_{(0^{10})}t_{10}^2t_9^{-2}t_8t_7^{2}t_6^{-3}t_5^4t_4^{-4}t_3^{-1}\}\}
\cup
\{C((a^{(01^9)}_{(0^{10})}t_{10}^2t_9^{-2}t_8t_7^{2}t_6^{-3}t_5^4t_4^{-4}t_3^{-1})^{\varepsilon})\}$.
For any reduced word $u=u_0t_4^{l_1}\cdots t_4^{l_n}u_n$, we have
$u\in (C(V'(Y'_1))$ if and only if
\begin{enumerate}
\item[(a)]\ $l_i<0\ (1\leq i\leq n),\ u_{i-1}$ is reduced on $Y'_1\backslash\{t_4\}$
and does not end with
$a^{(01^9)}_{(0^{10})}t_{10}^2t_9^{-2}t_8t_7^{2}t_6^{-3}t_5^4$,
\item[(b)]\ $l_i>0\ (1\leq i\leq n),\ u_{i-1}$ is reduced on $Y'_1\backslash\{t_4\}$
and does not end with $b_{(120^8)}^{\varepsilon},\
c_{(120^8)}^{\varepsilon}$,
\item[(c)]\ $u_n$ is reduced on $Y'_1$.
\end{enumerate}

Let $S$ consist of the above relations (4.1)-(4.20) and the trivial
relations in $H_4$. Let $X=((X_1^{\pm1}\dot{\cup}\{t_3^{\pm1},\
t_4^{\pm1},\
t_7^{\pm1}\})\dot{\cup}\{t_2^{\pm1}\})\dot{\cup}\{t_1^{\pm1}\}$ and
$$
t_i>t_i^{-1},\ t_1>t_2>t_3
>t_4>t_7>t_5>t_6>t_8>t_9>t_{10}.
$$
Then, with the first tower order on $X^*$, $S$ is an effective
standard Gr\"{o}bner-Shirshov basis. By Lemma \ref{l2.6}, $Irr(S)$
is an effective $k$-basis of the algebra $kH_4$. Since the canonical
forms of $H_4$ is $Irr(S)$, $H_4$ is a group with the effective
standard Gr\"{o}bner-Shirshov basis and the effective standard
normal form, and $H_1\leq H_2\leq H_3\leq H_4$ a tower of
HNN-extensions.
\end{example}

\section{Question}

In this section, we present a question related to the one-relator
groups.

\noindent{\bf Question}: Let $G_0\rightarrow
G_1\rightarrow\cdots\rightarrow G_n$ be a Novikov tower. For any
$G_i\ (0\leq i\leq n)$ and relations $A_ip=pB_i$ in $G_i$, assume
that the following conditions hold:
\begin{enumerate}
\item[(i)]\ $\{A_i\},\ \{B_i\}$ freely generate two subgroups of $G_{i-1}$
respectively;
\item[(ii)]\ $A_ip=pB_i$ can be presented as $\ A'_ixA''_ip=pB'_iyB''_i$,
where $x,\ y$ are two distinguishing letters of the highest weight
in the words $A_i$ and $B_i$ ($x$ is relative to $p$ and $y$ to
$p^{-1}$), respectively;
\item[(iii)]\ All the distinguishing letters relative to $p$ are different (also to $p^{-1}$ are
different).
\end{enumerate}
Can we  get an effective standard normal form for $G_n$?

We can see that all the one-relator groups mentioned in this paper
can be embedded into Novikov towers which satisfy the above
conditions and have the effective standard normal form.

\section{Regression to embedding into two-generator groups}

Let
\begin{eqnarray*}
G&=&gp\langle X| S\rangle,\\
G_1&=&gp\langle X,\ a,\ b| S,\
x_i=a^{-1}b^{-2i+1}a^{-1}b^{2i-1}ab^{-2i+1}ab^{2i-1}\rangle,
\end{eqnarray*}
where $X=\{x_i|1\leq i\leq n\}$ ($n$ may be infinite).

B. H. Neumann and H. Neumann  \cite{nn59} proved that $G$ can be
embedded into $G_1$.

Now, we  use the Magnus' method and Composition-Diamond lemma to
reprove the B. H. Neumann and H. Neumann's embedding theorem.

We may assume that $S$ is a Gr\"{o}bner-Shirshov basis for $G$.

Let $b=t,\ a_i=b^{-i}ab^i$ and rewrite the last defining relation as
$x_i=a_0^{-1}a_{-2i+1}^{-1}a_0a_{-2i+1}$. Then
$$
H_1=gp\langle X,\ t,\ a_j\ (-2n+1\leq j\leq 0)| S,\
x_i=a_0^{-1}a_{-2i+1}^{-1}a_0a_{-2i+1},\ a_jt=ta_{j-1}\ ( j-1\neq
0)\rangle
$$
and $G_1\cong H_1$ by $x\mapsto x,\ a\mapsto a_0,\ b\mapsto t$. Let
$a_{-2i+1}=t_i$. Then
$$
H_2=gp\langle X,\ t,\ t_i,\ a_{-2i+2}\ (1\leq i\leq n)| S,\
a_0t_{i}=t_{i}a_0x_i,\ a_{-2i+2}t=tt_{i},\ t_{i}t=ta_{-2i} (i\neq
n)\rangle
$$
and $H_1\cong H_2$.

Let $Y=\{X^{\pm1},\ t^{\pm1},\ t_i^{\pm1},\ a_{-2i+2}^{\pm1}\ |1\leq
i\leq n\}$. Define the tower order on $Y^*$ as Definition \ref{d2.2}
with $t>t^{-1}>t_i>t_i^{-1}>a_{-2j+2}>a_{-2j+2}^{-1}>x>x^{-1}$.

The relative standard relations in $H_2$ are as follows.
\begin{enumerate}
\item[(6.1)]\ $a_0^{\varepsilon}t_{i}=t_{i}(a_0x_i)^{\varepsilon},\
a_0^{-1}t_{i}^{-1}=x_it_{i}^{-1}a_0^{-1},\
 a_0x_it_{i}^{-1}=t_{i}^{-1}a_0$,

\item[(6.2)]\ $a_{-2i+2}^{\varepsilon}t=tt_{i}^{\varepsilon},\
t_{i}V(a_0x_i)t^{-1}=V'(a_0)t^{-1}a_{-2i+2},\
t_{i}^{-1}V'(a_0)t^{-1}=V(a_0x_i)t^{-1}a_{-2i+2}^{-1}$,

\item[(6.3)]\ $t_{i}V(a_0x_i)t=V'(a_0)ta_{-2i},\
t_{i}^{-1}V'(a_0)t=V(a_0x_i)ta_{-2i}^{-1},\
a_{-2i}^{\varepsilon}t^{-1}=t^{-1}t_{i}^{\varepsilon}$,
\end{enumerate}
where $V(a_0x_i)\leftrightarrow V'(a_0)$ by $a_0x_i\leftrightarrow
a_0$.

Let $R$ consist of the relations $(6.1)-(6.3),\ S$ and the trivial
relations of $H_2$. Since $S$ is a Gr\"{o}bner-Shirshov basis for
$G$ and there is no composition occurred among $(6.1)-(6.3)$ with
$S$, $R$ is a Gr\"{o}bner-Shirshov basis for $H_2$ with the tower
order. Therefore, by Composition-Diamond lemma, $Irr(S)\subset
Irr(R)$ and so, $G$ can be embedded into $H_2\cong G_1$.

G. Higman, B. H. Neumann and H. Neumann  \cite{hnn} also proved that
$G$ can be embedded into $G_2$, where
$$
G_2=gp\langle X,\ a,\ b| S,\
x_i=a^{-1}b^{-1}ab^{-i}ab^{-1}a^{-1}b^{i}a^{-1}bab^{-i}aba^{-1}b^{i}\rangle.
$$
Now, we reprove this result by using also the Magnus' method and
Composition-Diamond lemma.

 Let $a=t_1,\
b_j=a^{-j}ba^j$ and rewrite the last defining relation as
$x_i=b_{-1}^{-1}b_0^{-i}b_{1}^{-1}b_0^{i}b_{-1}b_0^{-i}b_{1}b_0^{i}$.
Then  we can get
$$
H_1=gp\langle X,\ t_1,\ b_j\ (-1\leq j\leq 0)| S,\
x_i=b_{-1}^{-1}b_0^{-i}b_{1}^{-1}b_0^{i}b_{-1}b_0^{-i}b_{1}b_0^{i},\
b_jt=tb_{j-1}\ ( j-1\neq 0)\rangle
$$
and $G_1\cong H_1$ by $x\mapsto x,\ a\mapsto t_1,\ b\mapsto b_0$.
Let $b_0=t_2,\ b_{(1j)}=b_0^{-j}b_1b_0^j,\
b_{(-1j)}=b_0^{-j}b_{-1}b_0^j$ and rewrite the defining relation as
$x_i=b_{(-10)}^{-1}b_{(1-i)}^{-1}b_{(-10)}b_{(1-i)}$. Let
$b_{(-10)}=t_3$. At last, we can get
\begin{eqnarray*}
H_2=gp\langle X,\ t_1,\ t_2,\ t_3\ b_{1-j}\ (1\leq j\leq n)&|& S,\
b_{(1-i)}^{-1}t_3=t_3x_ib_{(1-i)}^{-1},\
b_{(1(-i+1))}^{-1}t_2=t_2b_{(1-i)}^{-1},\\
&& b_{(10)}t_1=t_1t_2,\ t_2t_1=t_1t_3\rangle
\end{eqnarray*}
and $H_1\cong H_2$.

Let $Y=\{X^{\pm1},\ t_1^{\pm1},\ t_2^{\pm1},\ t_3^{\pm1},\
b_{1-j}^{\pm1}\ |1\leq j\leq n\}$. Define the tower order on $Y^*$
as Definition \ref{d2.2} with
$t_1>t_1^{-1}>t_2>t_2^{-1}>t_3>t_3^{-1}>b_{(1-j)}>b_{(1-j)}^{-1}>x>x^{-1}$.

The relative standard relations in $H_2$ are as follows.
\begin{enumerate}
\item[(6.4)]\ $b_{(1-i)}^{\varepsilon}t_3=t_3(x_ib_{(1-i)}^{-1})^{\varepsilon},\
b_{(1-i)}^{-1}t_3^{-1}=x_i^{-1}t_3^{-1}b_{(1-i)},\
 b_{(1-i)}x_i^{-1}t_3^{-1}=t_3^{-1}b_{(1-i)}^{-1}$,

\item[(6.5)]\ $b_{(1(-i+1))}^{\varepsilon}t_2=t_2b_{(1-i)}^{\varepsilon},\
b_{(1-i)}^{\varepsilon}t_2^{-1}=t_2^{-1}b_{(1(-i+1))}^{\varepsilon},\
b_{(10)}^{\varepsilon}t_1=t_1t_2^{\varepsilon}$,

\item[(6.6)]\ $t_{2}V(b_{(1-i)})t_1^{-1}=V'(b_{(1(-i+1))})t_1^{-1}b_{(10)},\
t_{2}^{-1}V'(b_{(1(-i+1))})t_1^{-1}=V(b_{(1-i)})t_1^{-1}b_{(10)}^{-1}
$,

\item[(6.7)]\ $t_2V(b_{(1-i)})t_1=V'(b_{(1(-i+1))})t_1t_3,\
 t_{2}^{-1}V'(b_{(1(-i)+1)})t_1=V(b_{(1-i)})t_1t_3^{-1},\
$,

\item[(6.8)]\ $
t_{3}^{-1}W(x_ib_{(1-i)}^{-1})t_1^{-1}=W'(b_{(1-i)}^{-1})t_1^{-1}t_{2)},\
t_3^{-1}W'(b_{(1-i)}^{-1})t_1^{-1}=W(x_ib_{(1-i)}^{-1})t_1^{-1}t_2^{-1}$,
\end{enumerate}
where $V(b_{(1-i)})\leftrightarrow V'(b_{(1(-i+1))})$ by
$b_{(1-i)}\leftrightarrow b_{(1(-i+1))}$,
$W(x_ib_{(1-i)}^{-1})\leftrightarrow W'(b_{(1-i)}^{-1})$ by
$x_ib_{(1-i)}^{-1}\leftrightarrow b_{(1-i)}^{-1}$.

Let $R$ consist of the relations $(6.4)-(6.8),\ S$ and the trivial
relations of $H_2$. It is clear that,  with the tower order, $R$ is
a Gr\"{o}bner-Shirshov basis for $H_2$. Therefore, $Irr(S)\subset
Irr(R)$ and so, $G$ can be embedded into $H_2\cong G_2$.

 \ \

\noindent{\bf Acknowledgement}: The authors would like to express
their deepest gratitude to Professor L. A. Bokut for his kind
guidance, useful discussions and enthusiastic encouragement.

 \ \

\end{document}